\documentclass[11pt]{article}
 \title{A two dimensional  adaptive nodes technique in irregular regions applied to  meshless-type methods}
\author{Kamal Shanazari \thanks{E-mail: kshanaz@liverpool.ac.uk.} \ and
        Mohammad Hosami  \thanks{E-mail: mohammad\_hosami@yahoo.com.} \\
         Department of Mathematics, \\
         University of Kurdistan, Sanandaj, IRAN.}
\date{ }
     \usepackage{amsmath}
     \usepackage{graphicx}
     \newtheorem{defi}{ Definition}
     \newtheorem{exam}{ Example    }
\newenvironment{@abssec}[1]{\vspace{.00in}\footnotesize \parindent .00in
         {\upshape\bfseries #1. }\ignorespaces } {\par\vspace{.00in}}
\newcommand\keywordsname{Key words}
\newcommand\AMSname{AMS subject classifications}
\newenvironment{keywords}{\begin{@abssec}{\keywordsname}}{\end{@abssec}}

\begin{document}  \maketitle
\begin{abstract}
We propose  a two dimensional (2D) adaptive nodes technique for irregular
regions. The method is based on equi-distribution principal and
dimension reduction. The mesh generation is carried out by first
producing some adaptive nodes in a rectangle based on
equi-distribution along the coordinate axes and then transforming
the generated nodes  to the physical domain. Since the produced mesh
is applied to the meshless-type methods, the connectivity of the
points is not used and only the grid points are important, though
the grid lines are utilized in the adapting process.
The performance  of the adaptive points is examined by considering a collocation
meshless method which is based on interpolation in terms of a set of  radial
basis functions.
 A generalized thin plate spline with sufficient  smoothness is used as a basis function and the arc-length is employed as a monitor in the equi-distribution process.
 Some experimental results will be presented to illustrate the
effectiveness of the proposed method. \\
\end{abstract}
\section{Introduction}
Mesh free methods are now well known in the numerical solution of
partial differential equations (PDEs). The most attractive feature of these techniques is that
  discretization  of the domain or boundary  is not required. This feature considerably reduces the computational complexity
 of the method. While the meshless techniques are often  divided into two categories: boundary (see for example \cite{wang}) and domain type methods, the current work concerns with the latter category. \\
 The use of radial basis functions (RBFs) for solving PDEs, first presented  by Kansa, is a fully
 mesh free approach and falls into
 the domain type matheds \cite{kansa_90, zhang_07}.
This method can be easily applied to the
case of higher dimensional spaces due to the nature of the RBFs.
Despite a good performance of RBFs in approximating  multi-variate functions,
they involve ill-conditioning,
especially for large scale problems. Another difficulty  concerns
their  computational efficiency, due to the dense
matrices arising from interpolation. To tackle the above difficulties some sort of
localization, such as domain decomposition methods (DDM)
\cite{dubal_94}  and compactly supported RBFs (CS-RBFs)
\cite{chen_02}, the most important of which was introduced by
Wendland \cite{wendland_95},  have been recommended.\\
In the DDM the domain is divided into some subdomains and the PDE
is solved for each subproblem followed by assembling the global
solution.  As a result, the ill-conditioning is avoided and the
computational efficiency is improved due to working with small
size matrices. On the other hand, using the CS-RBFs results in
sparse matrices, which again improve the conditioning and
computational
efficiency of the method.\\

This work involves a different approach which can still be used
with the above proposed methods. As was  highlighted before, in
using the classical RBFs, increasing the size of the problem
itself affects the conditioning. Consequently, reducing
the number of nodes can improve the conditioning. One way to
achieve this goal is to apply a set of adaptive nodes rather than
uniform ones. As is well known, the main idea in adaptive meshes is
to use a minimum number of nodes while still having the desired
accuracy.  This is achieved  by allocating more mesh points to
the areas where they are required. The adaptive mesh strategies
often fall into two categories: the equi-distribution principle
\cite{deboor_73}  and the variational principle \cite{winslow_67}.
The most popular technique, which has
been widely used in the literature, is based on the
equi-distribution strategy, which is also employed in this work. In
this approach  the mesh distribution is carried out in such a way
that some measure of error, called a monitor function, is
equalized over each subinterval. \\

 Much effort has been devoted to generating  adaptive meshes in two and three dimensional
spaces,  based on both equi-distribution and variational
principles (see for example \cite{olga, burg, rag}). In the
literature two major methods, namely  transformation \cite{chen2}
and dimension reduction \cite{sweby}, have been employed to
produce 2D meshes.
 The first category is based on mapping the physical domain into a
 simple domain with a uniform mesh  and it leads to solving a
 differential equation in order to obtain an adaptive mesh.
 In the dimension  reduction method, which is also employed in this paper, the equi-distribution process is reduced to a 1D case. \\

A method based on
equi-distributing along the grid lines in the coordinate
directions has been presented to produce mesh points in rectangular regions \cite{kamal}. Also a generalization of this method to
the case of three dimensions was proposed in \cite{kamal_08}.
However, due to the use of grid lines in the coordinate
directions, this method was limited to the case of rectangular and
cubic domains, respectively, in the case of 2D and 3D. The purpose
of the current work is to extend the above method to more general
cases with irregular boundaries in 2D. This is carried out by first,
generating some adaptive nodes in a rectangle, then mapping  the generated nodes  to the physical domain
employing a suitable transformation.
Of course, the mesh produced by the proposed method neither
precisely satisfies the equi-distributing condition nor concerns
about properties  such as orthogonality, which are often required
in numerical mesh-dependent methods. Instead, the mesh points
are suitable for any meshless
methods in  which the mesh points, rather than mesh lines, are important. \\

While a part of  researches  in the adaptive mesh community
concerns with constructing  monitor functions for different
applications \cite{huang}, the current study focus on the mesh
generation strategy using a well known monitor, namely arc-length
\cite{white}. In addition, we remark that the connectivity  of
the mesh are not used in this work, although they are used in the
mesh points
generation process. \\
 This paper is organized as follows. In section \ref{adapt2} the adaptive mesh technique in 1D is reviewed. A generalization of this method to the case of 2D is discussed in section \ref{adapt3}. The new mesh generation  method is presented in section \ref{adapt4}.
In section \ref{meshless} the collocation meshless method is
reviewed.
Some numerical results are given in section \ref{numer}.

\section{ Adaptive mesh } \label{adapt2}
We now introduce the concept of equi-distribution in the case of 1D.
\begin{defi}[Equi-distributing] \label{equi_d}
Let $M$ be a non-negative piecewise continuous function on $[a,b]$
and $c$ be a constant such that $ n=\frac{1}{c}\int _{a}^{b}M (x)
dx $ is an integer. The mesh
\[ \Pi : \quad a=x_{0}<x_{1}< \cdots <x_{n}=b, \]
is called equi-distributing (e.d.) on $[a,b]$ with respect to M
and $c$ if
\[ \int _{x_{j}}^{x_{j+1}} M(x) dx=c,  \qquad    j=0, 1, \cdots,
n-1.\]
\end{defi}

A  suitable algorithm to produce an e.d. mesh has been given in
\cite{kautsky}.  In  Definition \ref{equi_d} the function $f$,
often called a monitor, is  dependent on the solution of the
underlying PDE and its derivatives.
 The
arc-length monitor
\begin{equation} \label{arc}
 M= \sqrt{1 + u^2_x} \ ,
 \end{equation}
 which is used in this work
 has been widely used in the literature (see for example \cite{white,
  beckett_m}). The function $u$ in (\ref{arc}) is the solution  of the underlying PDE,  $x$ is the coordinate, in the direction of
  which the adaptivity is performed, and $u_x$ is the partial derivative with respect to
  $x$.  To find more details about the monitors for different applications see for
  example  \cite{carey_hum_81, huang}. \\
\section{Adaptive nodes in a rectangle} \label{adapt3}
\begin{figure}
\includegraphics{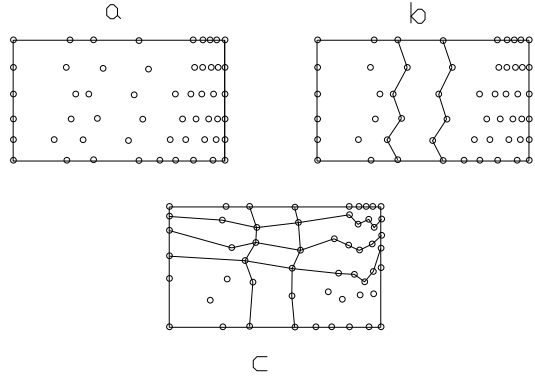} \vspace{5cm} \caption{The  three stages of the
adaptive mesh method are shown in Figures (a), (b) and (c)
respectively and in each direction the grid curves are displayed
with  the quadrilateral formed.}
 \label{3stg}
\end{figure}

A natural extension of Definition \ref{equi_d} to the case of 2D  is
as follows,
\begin{defi}[2D Equi-distributing] \label{equi_3d}
 Given a  2D  domain $\Omega$, a 2D adaptive mesh
 based on equi-distributing will be a mesh obtained by dividing
 the domain $\Omega$ into $n$ subdomains $\Omega_{i}$ such that
 $$\int\int_{\Omega_{i}} M(x,y)dxdy=constant, $$
where $M$ is a suitable monitor function.
\end{defi}
 Obviously an infinite
number of adaptive meshes based on Definition \ref{equi_3d} exist.
However, obtaining even one of these meshes can be a complicated
process. %
 A method  based on dimension reduction was proposed for a
 rectangular region in \cite{kamal}. Since the current work is a development
  of the above mentioned method, here it is briefly  reviewed.
 To do so,  we start with a uniform mesh in a rectangle in the
form
$$\{(x,y)|\ a_{1} \leq x \leq b_{1},\ a_{2}\leq y
\leq b_{2}\}.
$$
Let the underlying uniform mesh points be
$$
\{(x_{ij},y_{ij})|\ i=0,1,\cdots,N_{1},\ j=0,1\cdots,N_{2} \},
$$

where \begin{equation} \label{meshsize}
 x_{ij}=x_{i}=a_1 +
i\cdot h_{1},\ y_{ij}=y_{j}=a_2 + j\cdot h_{2},
\end{equation}
 and
$$h_{m}=\frac{b_{m}-a_{m}}{N_{m}}, \
\ \ \ m=1,2.$$%
 The equi-distribution  process is performed in three stages.  In the first stage, equi-distribution is performed
in the direction of the $x$-axis. More precisely, for each
horizontal line $y=y_{j}, \ \ j=0,1,\cdots,N_{2}$,\ we obtain the
new mesh
$$
(x^{'}_{ij},y_{ij})
$$
such that
$$
\int_{x^{'}_{ij}}^{x^{'}_{i+1j}} M_{x}(x,y_{ij})dx=constant,\ \ \
i=0,1,\cdots,N_{1}-1,
$$
where $M_{x}$ is the monitor function in the  $x$ coordinate
direction. Note that only the first coordinates of the points have
been changed (Fig. \ref{3stg}a). \\
In the second stage the equi-distribution is performed in the
vertical direction along the grid lines produced in the first
stage. Since the grid lines are curved,   the distribution is
performed  along the arc rather than the vertical coordinate.
Denoting  by $s$ the arc-length variable for each vertical grid
line $i=0,1,\cdots,N_{1},$ \ the distance $s_{ij}, \ \
j=0,1,\cdots,N_{2}$\ from \ $(x^{'}_{io},y_{io})$\ to \
$(x^{'}_{ij},y_{ij})$ \ along the vertical grid lines can be
evaluated piecewise linearly by ($ s_{io}=0$ )
\[ s_{ij}=s_{i(j-1)}+\|(x_{ij},y_{ij})-(x_{i(j-1)},y_{i(j-1)})\|. \]
Having the values of the monitor function
corresponding to \ $s_{ij}$, \ i.e. \ the value of \ $M_{y}$\ at
the points \ $(x^{'}_{ij},y_{ij}),$\ the new e.d.  mesh
$s^{'}_{i0}, s^{'}_{i1},\cdots,s^{'}_{iN_{2}}$\ is obtained by
$$
\int_{s^{'}_{ij}}^{s^{'}_{i(j+1)}} M_{y}(x,y)ds=constant, \ \
 i=0,1,\cdots,N_{1}.
$$
\begin{figure}
\includegraphics{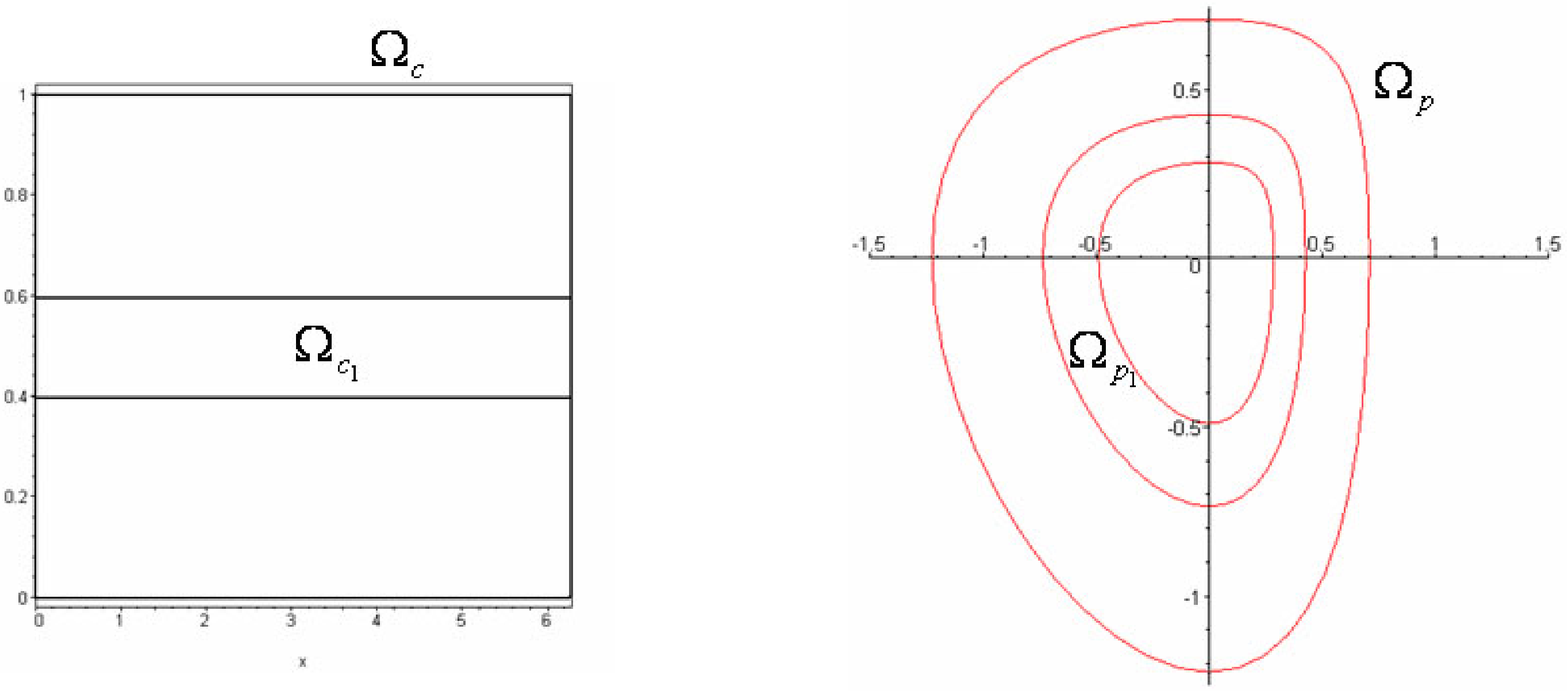} \vspace{8cm} \caption{The subregion $\Omega_{c_i}$
in $\Omega_c$ is mapped into $\Omega_{p_i}$ in $\Omega_p$.}
\label{omegacp}
\end{figure}
 The new values of  $s^{'}_{ij}$ can be used to
generate the new positions of the points on the grid lines since
the  piecewise linear representation of the underlying arc is
available (Fig. \ref{3stg}b). For more details see \cite{kamal}.
 A similar procedure is performed in the third stage along the
horizontal grid curves, with the monitor in the \ $x$\
coordinate direction again (Fig. \ref{3stg}c).\\
The mesh resulting from the above procedure  forms quadrilaterals
whose sides equi-distribute the
grid lines in the two coordinate directions (see Fig. \ref{3stg}c). \\
\section{Adaptive nodes in a non-rectangular domain} \label{adapt4}
In this section we present a new technique for generating
adaptive nodes in a simply connected domain, $\Omega_p$, bounded
by a closed curve $\Gamma$. For ease of explanation, we assume
that  $\Gamma$ can be given by the parametric equations
\[  \Gamma : x=g_1(\theta), \quad y=g_2(\theta), \quad 0\leq \ \theta  \leq 2 \pi , \]
where $\theta$ is the angle used in the polar coordinate system and measured in the conventional anti-clockwise
direction. The key to our proposed method is
to make use of a rectangular transformation and to produce an
adaptive mesh in the rectangle. We introduce a mapping
\begin{equation}  \label{trans} \psi : \Omega _c  \rightarrow \Omega _p
\end{equation}
 where
 \[\Omega_c = \{ (\theta,r), 0\leq  \theta  \leq 2 \pi, \  0 \leq r \leq 1 \}\]
  is a rectangle in the cartesian coordinate system $(\theta,r)$, referred to as computational domain,   $\Omega_p$ is the physical domain and the transformation is
defined by the following functions
\begin{equation}  \label{transform}   \left  \{ \begin{array}{l}  x=r g_1(\theta),
  \\ y= r g_2(\theta).
  \end{array}  \right.  \end{equation}
 One can
easily see that the above transformation corresponds  the line $r=1$ in
$\Omega_c$ to the boundary $\Gamma$ in $\Omega_p$. In addition, any
rectangular subregion $\Omega_{c_i}$ in $\Omega_c$ is mapped into a
subregion $\Omega_{p_i}$ in $\Omega_p$ as shown in Fig. \ref{omegacp}. \\
 The method described
in  section \ref{adapt3} and used to produce grid points in a rectangle, is
now applicable to the computational domain $\Omega_c$. Before
explaining the adapting  technique, we propose a method to construct
a uniform mesh in $\Omega_p$ by the use of the transformation in (\ref{trans}).

We start with  a uniform mesh in $\Omega_c$ in the cartesian
coordinate system $\theta r$ as follows,
\[ \theta_0 < \theta_1 < \dots < \theta_{n-1} = 2 \pi -h_\theta,\]
\[  h_r = r_1<r_2< \dots < t_{m-1} < r_m =1,\]
where
\[ h_{\theta} = \frac{2 \pi}{n}, \ \theta_j =j h_\theta,  \ j=0,1, \dots, n-1,\]
and
\[ h_r=\frac{1}{m}, \ r_i= i h_r, \ i=1,2, \dots, m.\]
\begin{figure}
\includegraphics{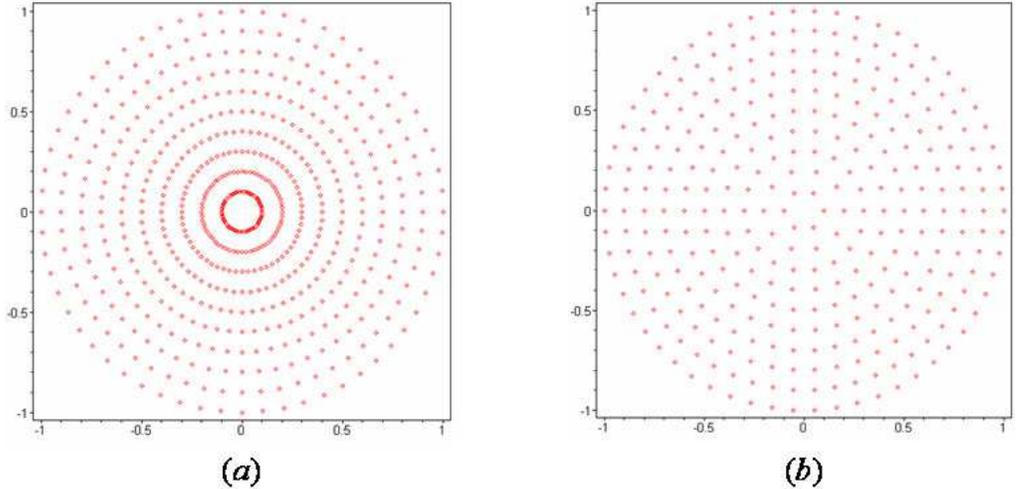}
\vspace{6cm}
\caption{The mesh points transformed from $\Omega_c$ into a circle before and after refinement
are displayed in Figures (a) and (b) respectively.}
\label{clus}
\end{figure}
  Using equations (\ref{transform})
the mesh points $(\theta_j, r_i)$ are transformed to the points
$(x_i,y_j)$ in the $xy$ coordinate system, as depicted in Fig. \ref{clus}.
Unlike the mesh in $\Omega_c$, the new mesh in $\Omega_p$  is non-uniform, since each line $r=r_i$ in
$\Omega_c$ corresponds to a closed curve in $\Omega_p$ with the same
number of nodes, n. Therefore, the mesh will be clustered around the
center, especially, for small values of $r$. This feature
affects the quality of the grid points.
 To avoid this difficulty, below  we
suggest a refinement process to  construct a roughly uniform mesh referred to as a uniform mesh in this paper.
Although there might be some easier way to produce a set of uniform mesh points
 in an irregular region, the following method
will be utilized in the adapting technique later.
\subsection{A uniform mesh in $\Omega_p$} \label{uniform}
\begin{figure}
\includegraphics{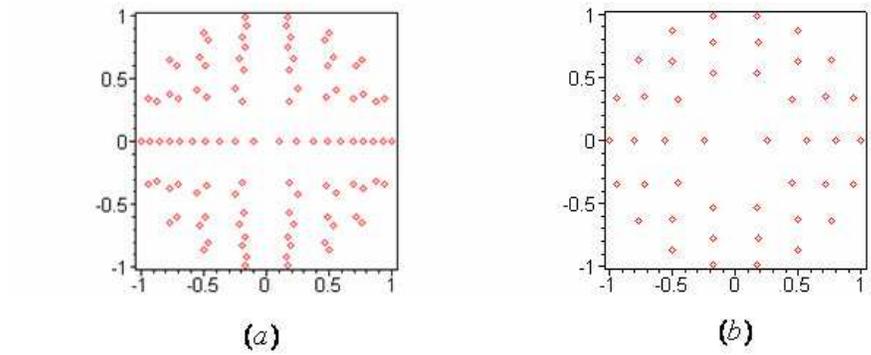}
\vspace{5cm}
\caption{The  points produced by an arbitrarily selected value of $n_r$ in (a) and  a suitably
selected value of $n_r$ based on relation (\ref{nr}) in (b) are displayed.}
\label{refinr}
\end{figure}

In order to obtain a uniform mesh in $\Omega_p$, we now refine
the above mesh points obtained by the transformation. The
refinement process is performed by modifying the number of nodes
on each closed curve, $r=r_i$ based on a suitable criteria. For
instance, we suggest this number of points to be selected in such
away that the distance between the adjacent points on each closed
curve be the
same as that on the boundary.  we propose  the following process which is based on the perimeter of the  closed curves.\\
 The perimeter of the $i$th curve can be approximated by the perimeter of a circle whose radius is evaluated by the average length of the position vectors of the boundary points, i.e.

\[ R_i= \frac{1}{n} \sum_{j=0}^{n-1} R_{ij}, \ i=1, \dots, m,  \] where

\[ R_{ij}= \sqrt{x^2_{ij} +y^2_{ij}}, \ \ x_{ij}=r_ig_1(\theta_j), \ \ y_{ij}=r_ig_2(\theta_j). \]
The approximate number of nodes for the $i$th curve can be therefore obtained by
\begin{equation}  \label{nthetai} n_{\theta_i} = \frac{p_i}{\Delta s}, \end{equation}
 where $p_i= 2 \pi R_i$ is the perimeter of the circle and $\Delta s$ is the minimum distance between the adjacent boundary nodes.
 Another difficulty may arise due to the value of $n_r$. In fact, if the values of $n_\theta$ and $n_r$ are arbitrarily selected, the points may be clustered along the lines $\theta=\theta_j$ (see Fig. \ref{refinr}). To treat this issue
 we propose the value of $n_r$ to be selected based on $n_\theta$ using a criteria similar to the above mentioned process. For instance, for a given value of $n_\theta$, $n_r$ can be determined  such that
 \begin{equation} \label{nr}
 \frac{p}{R} = \frac{n_\theta}{n_r},
 \end{equation}
 where $p$ is the approximate perimeter of the boundary and $R$ is the average lengthes of the position vectors.
 A  clustered  and a refined mesh obtained by the above process
 are displayed in  Figures  \ref{clus}a and  \ref{clus}b, respectively.\\
The suitably selected $n_{\theta_i}$  and $n_r$  will be also used in  adapting  mesh points later.
\subsection{Adaptive nodes in  non-rectangular regions}
\begin{figure}
\includegraphics{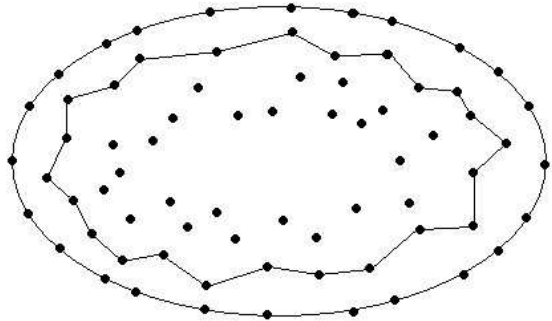}
\vspace{6cm}
\caption{The grid curves in $\Omega_p$ resulted from the transformation of the grid lines
 in $\Omega_c$ in the third stage of adapting nodes.}
\label{glin}
\end{figure}

We now present a method to produce adaptive nodes in the physical domain $\Omega_p$.
 The key to this new approach is to make use of the
adaptive mesh technique described in section \ref{adapt3} and the
transformation (\ref{trans}). More precisely, some adaptive points are first
prepared in the computational domain $\Omega_c$ and then
transformed into the physical domain $\Omega_p$. Fig. \ref{glin} illustrates
how, for instance, the third stage of  adapting mesh in $\Omega_c$ is related to
$\Omega_p$, i.e. equi-distributing along the grid lines for the coordinate $\theta$. \\

 Although, the mesh points produced by the above
method is somehow adaptive, the concentration of the points
around   the center is an disadvantage as previously discussed.
To overcome this
difficulty, below we propose a refining mechanism
in performing the 3-stage adaptive algorithm. \\
We start with a uniform mesh in $\Omega_c$ and   perform the the
first two stages of the adapting method. More precisely, the equi-distribution
in the horizontal direction and the vertical grid lines are
accomplished (Fig. \ref{3stg}a and  \ref{3stg}b). But, the final stage is done in a different manner.
We suggest a combination of the third stage of the adapting
technique with the refinement  process discussed in section \ref{uniform} to
avoid clustering the mesh points. As noted before, the
density of the mesh around the center can be avoided by modifying
the number of nodes along the closed curves. Recalling the
refinement process to obtain a uniform mesh, a suitable
number of points for the $i$th closed curve, $n_{\theta_i}$ was
 computed in (\ref{nthetai}). We can either use the same number
of points or a new value of $n_{\theta_i}$  based on the new
perimeter of the curves resulted from the first two stages. \\
Having obtained the horizontal grid lines by the two-stage procedure
and found the $n_{\theta_i}$ for each curve, we can equi-distribute
  $n_{\theta_i}$ points
  along the horizontal grid curves used in the third stage of adapting method in Fig. \ref{3stg}c.
  \begin{figure}
\includegraphics{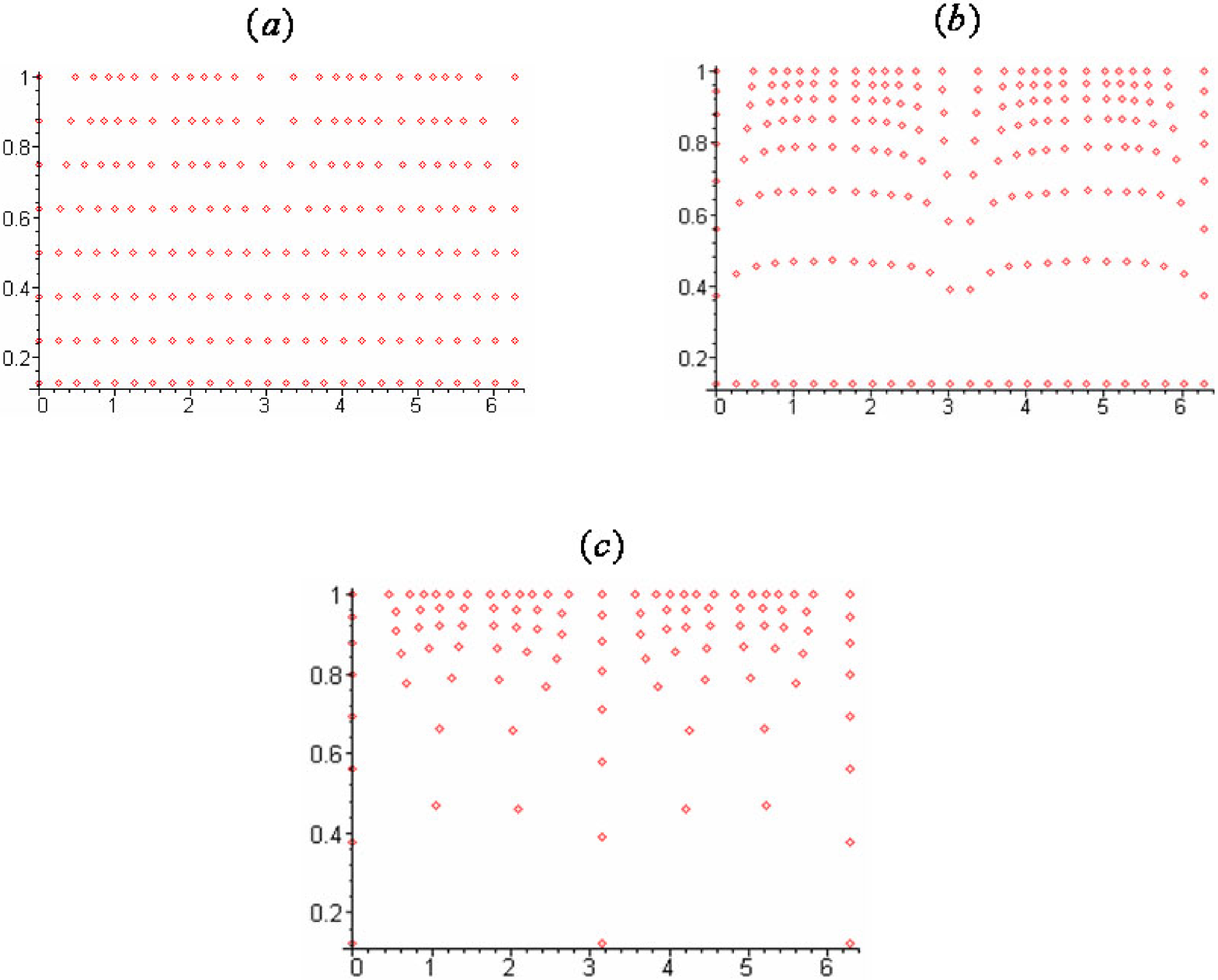}
\vspace{11cm}
\caption{The illustration of the 3-stage adapting mesh in $\Omega_c$ with a refined number
 of nodes along each grid cure in the horizontal direction in (c) is displayed.}
\label{3exp}
\end{figure}
  Fig. \ref{3exp} illustrates the 3-stage adapting nodes in the rectangular domain. It is observed
  that in the  third stage (Fig. \ref{3exp}c) a refined number of nodes, $n_{\theta_i}$ are distributed along the
  horizontal grid curves.
\section{Meshless methods} \label{meshless}
In this section we first
introduce the  RBFs and then describe their  application to the
numerical  solution of  PDEs  based on  collocation meshless method.
\subsection{Radial basis functions}
RBFs are known as the natural extensions of  splines to
multi-variate interpolation. Suppose the  set of points
$$
\{x_{i}\in\Omega|i=1,2,\cdots,N\},
$$ is given,
where $\Omega$\ is a bounded domain in $R^n$.
 The radial
function\ $\varphi:\Omega\longrightarrow R$\ is used to construct
the approximate function
$$
s(x)=\sum_{k=1}^{N}\alpha_{k}\varphi(\|x-x_{k}\|),
$$
which interpolates an unknown function $f$ whose values at
$\{x_{i}\}_{i=1}^{N}$ are known.\ $ \| . \| $ represents the
Euclidean norm. The unknown coefficients $\alpha_{k}$\ are
determined such that the following \ $N$\ interpolation conditions
are satisfied,
$$
f(x_{i})=s(x_{i})=\sum_{k=1}^{N}\alpha_{k}\varphi(\|x_{i}-x_{k}\|),
\ \ \ i=1,\cdots,N.
$$
There is a large class of interpolating RBFs \cite{powell} that
can be used in  meshless methods. These include the linear $1+r$,
the polynomial $P_k(r)$, the thin plate spline (TPS) $r^ 2 \log
r$,
 the Gaussian $\exp (-r^2 / \beta^2)$, and the
multi-quadrics $\sqrt{\beta^2 +r^2}$ (with $\beta$  a constant
parameter). In this paper we  employ  a generalized TPS, i.e.  $r^4 \log r$
which is a  particular  case of $r^{2k} \log r$  (k=2).
These RBFs, augmented by some polynomials, are known as the natural
extension of the  cubic splines to the case of 2D. In fact, they are obtained by minimizing a $H^m$
semi-norm over all interpolants for which the semi-norm exists. The
theoretical discussion of these RBFs has been presented in
\cite{duchon_77}.
\subsection{ Collocation meshless method} \label{meshfree}
We describe the collocation method for a general case of PDEs in
the form
\begin{eqnarray}
\label{eq1} Lu=F,
\end{eqnarray}
where \ $L=[L_{1},\cdots,L_{N}]^{T}$\ represents a vector of
linear operations and   $F=[f_{1},\cdots,f_{N}]^{T}$\ denotes a
vector containing the right hand sides of the equations.\ For
instance,\  Poisson's equation with a Dirichlet boundary condition
$$
\begin{array}{cc}
  \Delta u=f, & \text{in} \ \Omega, \\
  u=g, &  \text{on} \ \partial\Omega, \\
\end{array}
$$
is a very simple case of equation (\ref{eq1}) where \
$L=[\Delta,I]^{T}$,  \
$F=[f,g]^{T}$ and the operators $\Delta$ and $I$ act on the domain $\Omega$ and the boundary $\partial \Omega$ respectively. \\
The collocation method is simply to express the unknown function\
$u$\ in terms of the\ $RBFs$\ as
\begin{eqnarray}
u(x)=\sum_{k=1}^{N}\alpha_{k}\varphi(\|x-x_{k}\|), \label{eq2}
\end{eqnarray}
and determine  the unknowns \ $\alpha_{k}$\  in such a way that
(\ref{eq2}) satisfies equation (\ref{eq1}) for all interpolation
points. Substituting (\ref{eq2}) in equation (\ref{eq1})\ and
imposing the \ $N$\ essential conditions of the collocation method
lead to a linear system of equations whose coefficient matrix
consists of \ $N$\ row blocks, the entries of which are of the
form
$$
A^{\mu}_{ij}=L_{\mu}\varphi(\|x-x_{j}\|)|_{x=x_{i}},\ \ \
i=1,\cdots,N_{\mu},\ j=1,\cdots,N,
$$
where $N_{\mu}$ indicates the number of nodes associated with the
operator $L_{\mu}$.\\
The above collocation method is referred to as a non-symmetric
collocation method due to the non-symmetric coefficient matrices.
The invertibility  of the coefficient matrix can not be guaranteed
\cite{fassh_00},\ although in most cases a non-singular matrix is
expected. To tackle this difficulty the symmetric collocation
method,\ motivated by Hermitian interpolation, has been suggested.
This method leads to  symmetric matrices and proves the
non-singularity of interpolation matrices at the expense of double
acting the operators on the $RBFs$ \cite{fassh_97}. Since this work
is not concerned with the singularity of the matrices the
non-symmetric case will be implemented. Of course, the main
technique  discussed in this paper can be applied to the other case
as well.
\section{Numerical Results} \label{numer}
We now examine the effectiveness of the  mesh generation technique  by applying the collocation meshless method
to some PDEs. In each case, the PDE is solved
 with some equally spaced nodes (roughly uniform)
and  adaptive nodes generated by the new method and the results are
 compared.
 As previously noted, $\phi(r)= r^4 \log r$, is used as a basis function.
In each example, $M$ test points, which do not coincide with the
interpolation nodes, are randomly selected and a  root mean square
(RMS) error at these points is evaluated by
\[
  \text{RMS error} =
\sqrt{\sum_{i=1}^M (u_{apr,i}- u_{ex,i})^2/M}\]
  where $u_{apr,i}$
and $u_{ex,i}$ denote the approximate and exact values of $u$,
respectively, at a test point $i$. \\

Equations (\ref{transform}) are used to evaluate the  arc-length monitors
$$M_\theta = \sqrt{1+u_\theta^2} \ \ \ \text{and} \ \ \ M_r= \sqrt{1+u_r^2}, $$
 respectively for equi-distributing in the  $\theta$ and $r$ coordinate  directions where
$$u_\theta=\frac{\partial u}{\partial x} \frac{dx}{d \theta} + \frac{\partial u}{\partial y} \frac{dy} {d \theta}, \ \ \ u_r=\frac{\partial u}{\partial x} \frac{dx}{d r}+ \frac{\partial u}{\partial y} \frac{dy}{dr}$$.

 We consider the Poisson equation
\begin{equation}\label{poissoneq}
\begin{array}{ll}
  \Delta u= f(x,y) & \text{in} \ \Omega \\
   u(x,y)=g(x,y) , & \text{on} \ \partial \Omega, \\
\end{array} \end{equation}
for  different cases with various solutions and regions.
\begin{exam} \label{ex1}
We solve equations (\ref{poissoneq}) in the case of
\[ \begin{array}{l}
f(x,y)=4e^{x^2+y^2}+4(x^2+y^2)e^{x^2+y^2}
\\ \\
g(x,y)=e^{x^2+y^2}
\end{array}  \]
 and $\partial \Omega$ is an ellipse whose equation is $x^2 + 4 y^2 -1=0$. \\
In this case the exact solution is given by  $u(x,y)= e^{x^2 +y^2}$.
\end{exam}
\begin{table} \caption{Error values for Example \ref{ex1}}
\begin{center}
\label{ext1}
\begin{tabular}{ccccccccc}  \hline\hline
$n_{\theta}, \ n_r$ & 25,4 & 31,5 & 37,6& 43,7 & 50, 8&56,9 & 63,10& 69,11 \\ \hline
  $n$ &60 & 86& 124 & 162 &212 &266&328&398 \\ \hline
  Adaptive & 1.04E-3  &  8.47E-5  & 5.15E-5 & 2.31E-6& 4.81E-6& 5.14E-7& 4.75E-7 &7.67E-8  \\
 Uniform & 2.53E-3 & 9.24E-4 & 3.42E-4 & 1.36E-4 & 5.57E-5 & 2.37E-5 & 1.09E-5 &4.80E-6
\end{tabular}
\end{center}
\end{table}
\begin{exam}
\label{ex2}
\[ \begin{array}{l}

f(x,y)=-4e^{(4-x^2-y^2)}+4(x^2+y^2)e^{(4-x^2-y^2)}
\\ \\
g(x,y)=e^{(4-x^2-y^2)}
\end{array} \] and $\partial \Omega$ is the same as that in Example \ref{ex1}. \\
The exact solution is given by  $u(x,y)= e^{(4-x^2 -y^2)}$.
\end{exam}
\begin{table} \caption{Error values for Example \ref{ex2}}
\begin{center}
\label{ext2}
\begin{tabular}{ccccccccc}  \hline\hline
$n_{\theta}, \ n_r$ & 25,4 & 31,5 & 37,6& 43,7 & 50, 8&56,9 & 63,10& 69,11 \\ \hline
  $n$ &60 & 86& 124 & 162 &212 &266&328&398                           \\ \hline
  Adaptive & 3.055E-4  & 3.76E-4  & 1.42-4 & 6.32E-5& 1.91E-5& 7.23E-6& 2.50E-6 &8.43E-7                        \\
 Uniform &2.89E-4 & 3.75E-4 & 2.26E-4 & 1.03E-4 & 4.32E-5 & 1.78E-5 & 7.83E-6 &3.13E-6
\end{tabular}
\end{center}
\end{table}
\begin{exam} \label{ex3}
\[ \begin{array}{l}
f(x,y)=4 \\
\\
g(x,y)=x^2+y^2
\end{array}  \] and the boundary is given by the parametric equations:
\[ \left \{ \begin{array}{l}
x=cos \theta \sqrt{1-cos \theta /2} \\
\\
y=sin \theta \sqrt{1-sin \theta/2}\\
\end{array}  \right. \]
The exact solution is given by  $u(x,y)= x^2 +y^2$.
\end{exam}
\begin{table} \caption{Error values for Example \ref{ex3}}
\begin{center}
\label{ext3}
\begin{tabular}{ccccccc}   \hline\hline
$n_{\theta}, \ n_r$ & 25,4 & 30,5 & 45,8& 55,9 & 65,11&75,12 \\ \hline
  $n$ &60 & 78& 142 & 262 &366 &460                           \\ \hline
  Adaptive & 1.13E-3  &  3.80E-5  & 2.54E-5 & 3.21E-6& 3.80E-7& 1.81E-7   \\
 Uniform & 1.85E-3 & 7.62E-4 & 5.40E-5 & 2.27E-5 & 5.34E-6 & 2.46E-6
\end{tabular}
\end{center}
\end{table}
\begin{figure}
\includegraphics{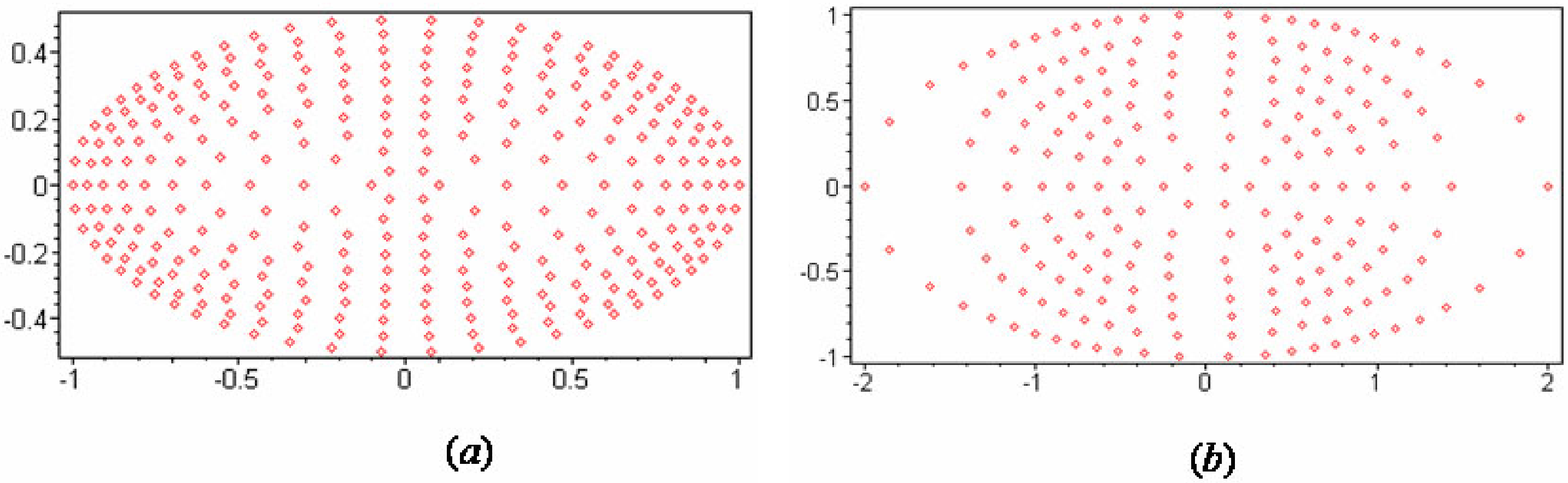} \vspace{9cm} \caption{The adaptive nodes generated
for the solutions of Examples \ref{ex1} and \ref{ex2} are
displayed in (a) and (b) respectively.} \label{exf12}
\end{figure}
\begin{figure}
\includegraphics{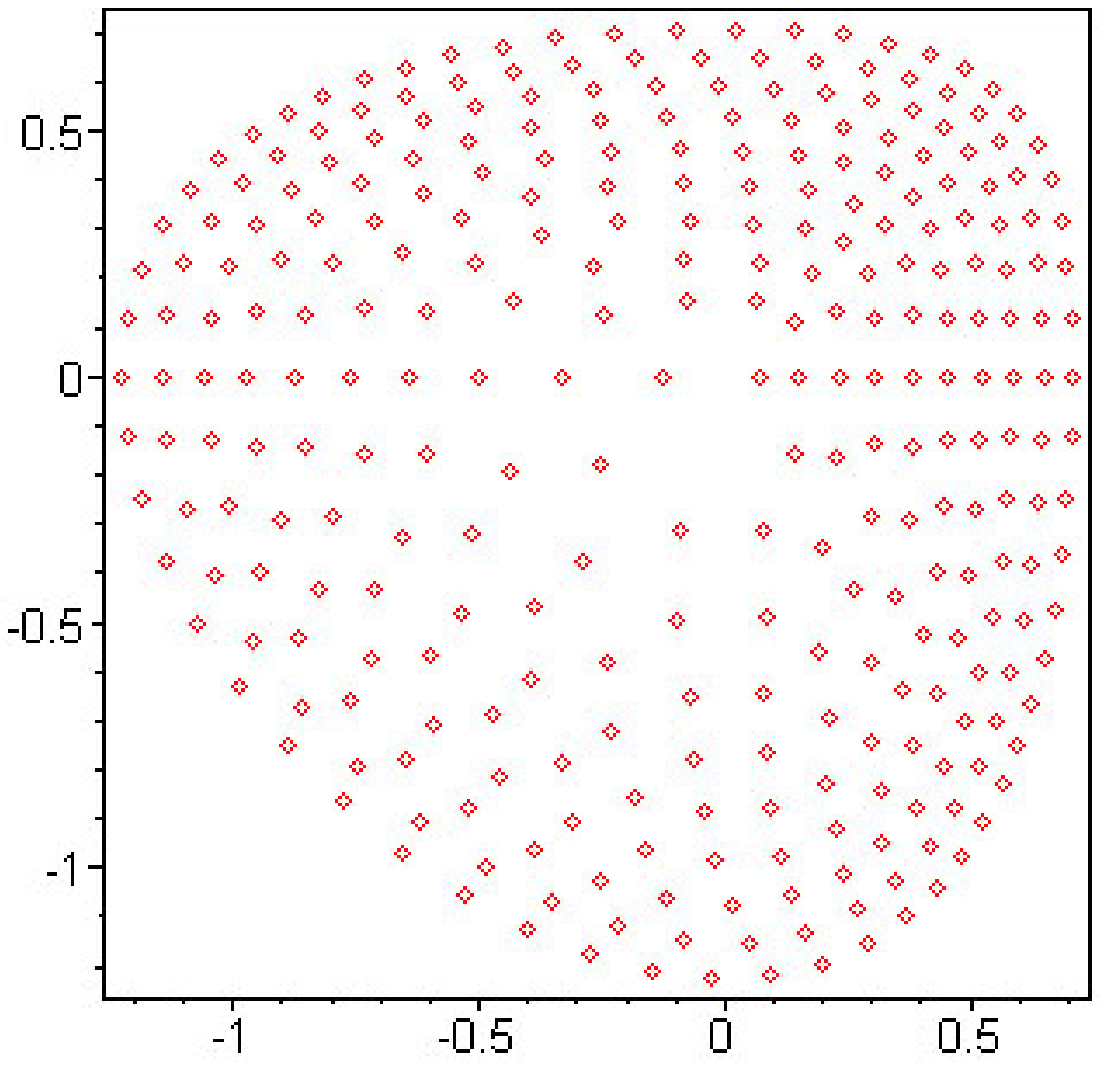}
\vspace{8cm}
\caption{The adaptive nodes generated for the solution of Example \ref{ex3}.}
\label{exf3}
\end{figure}
The adaptive nodes generated for the above PDEs are displayed in  Figures \ref{exf12}a, \ref{exf12}b and \ref{exf3}, respectively. One can expect, for the solution of Example \ref{ex1}, the mesh points to be more dense close to the boundary and, for the solution of Example \ref{ex2}, the mesh points to be more dense around the center. These are exactly what we observe in the figures. Therefore, the adaptive nodes are in a good agreement with the solutions. A similar situation is observed in Fig. \ref{exf3} for the solution of Example \ref{ex3}.\\

The numerical errors for the above three examples are listed, respectively, in Tables \ref{ext1}, \ref{ext2} and \ref{ext3}, where $n_\theta$ is the number of boundary points, $n$ is the total number of collocation points and $n_r$ is the appropriate number of divisions for the coordinate $r$  based on relation (\ref{nr}).\\
In all the examples, the numerical results demonstrate a considerable reduction in the error in the case of using the adapting nodes produced by the new adaptive method. In particular, Table \ref{ext1} shows that  the error in the case of using $398$ uniform nodes is nearly the same as that in the case of using $212$ adaptive nodes. Moreover, the error in the case of using $398$ uniform nodes is $63$ times larger than that in the case of using the same number of adaptive mesh points.
\section{Conclusion}
A new  adaptive mesh technique was presented  for  irregular regions in  two  dimensional cases based on equi-distribution.
The method was based on first, generating some adaptive nodes in a rectangle  and then mapping the produced nodes to the physical domain using a suitable transformation. \\

The new method was examined by considering some PDEs solved by a collocation meshless method and the
results demonstrated  considerable reduction in the
error values. Since the current work was not involved with
constructing a monitor for the underlying method, a general
monitor function, arc-length, was employed to produce the adaptive
nodes. Of course, using a suitable
monitor for the underlying method could have resulted in more
improvement in the numerical results \\

The proposed method was implemented   for 2D regions whose
boundary was given by  parametric equations in terms of polar
coordinate system. For a more general case of irregular
boundaries, we suggest a piecewise interpolation, in terms of
some scattered points on the boundary, to make a  set of
parametric equations. Moreover, an extension of the new method to
the case of 3D is currently under study and will be reported in a
near future. 
\end{document}